\documentclass[12pt,english]{article}
\usepackage[T1]{fontenc}
\usepackage[latin9]{inputenc}
\usepackage{textcomp}

 

\usepackage{babel}

\begin{document}

\title{New special functions solving nonlinear autonomous dynamical systems}

\date{}

\maketitle
\vspace{-7mm}

\hrule

\begin{center}
{\large L\'eon Brenig}%
\footnote{{\footnotesize Departement of Theoretical and Mathematical Physics,
Faculty of Sciences, Université Libre de Bruxelles (ULB), CP.231,
1050 Brussels, Belgium, email : lbrenig@ulb.ac.be}%
}{\large{} } 
\par\end{center}


\begin{abstract}
A general solution is found for a large class of time continuous autonomous
nonlinear dynamical systems, the so-called quasi-polynomial systems.
This solution is expressed in terms of a new type of special functions defined via their Taylor series. The coefficients of these Taylor
series are generated by a tensor that generalizes
the factorial function and has a  combinatorial meaning.
The existence of these functions raises the
question of the relation between them and the chaotic behaviour of
the solutions that may appear for the quasi-polynomial dynamical systems.
\end{abstract}

\section*{Introduction}

This article is intended to draw the attention of specialists in the fields of special functions, combinatorics and nonlinear dynamics on a new class of special functions that solve a broad class of nonlinear dynamical systems.\\ \\
Let us define an autonomous dynamical system in continuous time as
a system of coupled nonlinear ODEs of the form \begin{equation}
\dot{x}_{i}=f_{i}(x_{1},\cdots,x_{n})\label{eq1}\end{equation}
 where the dot denotes the time derivative, the index $i$ runs from
$1$ to $n$, the dependent variables $x_{i}$ are real, the functions
$f_{i}$ are enough regular to ensure the existence and unicity of
solutions, and the boundary conditions are given at the initial time.\\

For such systems there does not exist a universal structure of the
general solution. This is in strong contrast with the particular case
of linear autonomous dynamical systems for which the general solution
of the Cauchy problem is given by the exponential of the constant
matrix that characterizes each of these systems.

One reason for this lack of universality is the infinite diversity
of possible functional forms of the functions $f_{i}$ in equations
(\ref{eq1}). Moreover, the nonlinearity of the equations entails
a factorial explosion in the coefficients of the Taylor series of
the solutions: A coefficient of order $k$ involves a sum of $k!$
terms that depend on the specific form of the functions $f_{i}$.
The recursion relations between these coefficients are as difficult
to solve as the original ODEs. Hence, one should not wonder that in
most cases it is impossible to find a closed form structure for the
general $k$ -th order Taylor coefficient of such solutions.\\

In this article, we show that for the class of quasi-polynomial systems
defined in the next chapter, a universal form of the general solution
to the Cauchy problem can be found. Moreover, as has been shown by
E.Kerner \cite{7}, most systems of non-polynomial nonlinear ODEs
relevant for physics, chemical kinetics and generally for mathematical
modeling in natural and social sciences can be reshaped in the form
of systems of polynomial ODEs. The latter, in turn, are a sub-class
of the quasi-polynomial systems of ODEs. This general property, thus,
extends our general solution to most nonlinear systems of ODEs defining
dynamical systems.

\section*{The quasi-polynomial dynamical systems }

Let us define the sub-class of dynamical systems (\ref{eq1}) for
which the functions $f_{i}$ are quasi-polynomials, i.e. finite sums
of monomials involving powers of the dependent variables that can
be real numbers. Such monomials we shall call quasi-monomials. \\The
systems of this class can be cast into a useful standard notation,
the so-called quasi-polynomial (QP) representation\cite{1,2}: \begin{equation}
\dot{x}_{i}=x_{i}\sum_{j=1}^{N}A_{ij}\prod_{k=1}^{n}x_{k}^{B_{jk}}\qquad\mbox{ for }i=1,\ldots,n\label{eq2}\end{equation}
 where $N$ refers to the number of quasi-monomials in the set of
variables $x_{j}$ in the right-hand side, and $A$ and $B$ are constant
rectangular matrices respectively $n\times N$ and $N\times n$ with
entries that are real numbers. The presence of the factor $x_{i}$
in front of the right-hand side of the above equation is essential
and constrains the definition of both matrices. However, it does not
restrict the generality of the class of systems (\ref{eq2}) as any
polynomial or quasi-polynomial can be represented in that form. The
presence of the factor $x_{i}$ underlines the importance of the logarithmic
time derivative of the functions $x_{i}(t)$ which plays a fundamental
role in the present approach. \\

Other authors have found independently the same form (\ref{eq2})
and derived from it theoretical results in the fields of chemical
reactions and ecological systems\cite{3,4,5}. These works were mainly
concerned with models describing complex networks of interacting entities
and they were essentially devoted to the study of the stability properties
of some of their solutions. More generally, the QP class covers almost
all the systems of interest for mathematical modeling in natural and
social sciences. \\

A fundamental feature of the above QP representation of a dynamical
system is its covariance under the group of quasi-monomial transformations.
These are transformations of the dependent variables defined as follows:
\begin{equation}
x_{i}=\prod_{k=1}^{n}\tilde{x}_{k}^{C_{ik}}\qquad\mbox{ for }i=1,\ldots,n\label{eq3}\end{equation}
 where the matrix $C$ is any invertible, constant, real, square $n\times n$
matrix. In order to avoid unrelevant discussions about the existence
of the inverses of these transformations and their differentiability
we limit our scope to systems and to initial conditions such that
their solutions remain in the positive cone. However, it can be proven
that most of the following results continue to be valid for systems
not fulfilling this restriction. \\

It can be easily shown that under transformations (\ref{eq3}) the
system (\ref{eq2}) becomes: \begin{equation}
\dot{\tilde{x}}_{i}=\tilde{x}_{i}\sum_{j=1}^{N}\tilde{A}_{ij}\prod_{k=1}^{n}\tilde{x}_{k}^{\tilde{B}_{jk}}\qquad\mbox{ for }i=1,\ldots,n\label{eq4}\end{equation}
 with the following rule of transformations for the matrices A and
B: \begin{equation}
\tilde{A}=C^{-1}A\label{eq5}\end{equation}
 and \begin{equation}
\tilde{B}=BC\label{eq6}\end{equation}
 where the products are matrix products. The identity of form of equations
(\ref{eq2}) and (\ref{eq4}) clearly exhibits the covariance of the
QP equations under the quasi-monomial group of transformations. Let
us stress that: \begin{equation}
\tilde{B}\tilde{A}=BA\label{eq7}\end{equation}
 in other words, the $N\times N$ matrix $BA$ is an invariant of
the quasi-monomial transformations (\ref{eq3}). This means that the
whole set of QP-systems is divided into equivalence classes labeled
by such $N\times N$ matrices. The fundamental matrix $BA$ is related
to the existence of a canonical form in each equivalence class as
we now show.\\

Indeed, under a particular quasi-monomial transformation, any QP system
can be brought to a canonical form, the so-called Lotka-Volterra system
of ODEs \cite{1,2,6}: \begin{equation}
\dot{x}_{i}=x_{i}\sum_{j=1}^{N}M_{ij}x_{j}\qquad\mbox{ for }i=1,\ldots,N\label{eq8}\end{equation}
 where $N$ is the number of monomials of the original QP system (\ref{eq2})
and the square $N\times N$ matrix $M$ is equal to the matricial
product $BA$ of the two matrices $B$ and $A$. There also exists
a second canonical form dual of the above one, which until now has
been much less studied \cite{1} but is also characterized by the
same matrix $BA$. \\

In a given equivalence class, the passage to the Lotka-Volterra form
is most easily shown in the particular case of QP systems for which
the number $N$ of quasi-monomials is equal to the dimension $n$
of the system (\ref{eq2}), the so-called square QP-systems. For the
square QP-systems the matrices $A$ and $B$ are of course square
matrices. In this case, and if the matrix $B$ is not singular, the
quasi-monomial transformation (\ref{eq3}) in which the matrix $C$
is chosen equal to $B^{-1}$ leads to the equation (\ref{eq4}) with
\begin{equation}
\tilde{A}=BA\label{eq9}\end{equation}
 and \begin{equation}
\tilde{B}=I\label{eq10}\end{equation}
 where $I$ is the identity matrix. Thus, the transformed equation
is now: \begin{equation}
\dot{x}_{i}=x_{i}\sum_{j=1}^{N}{(BA)}_{ij}x_{j}\qquad\mbox{ for }i=1,\ldots,N\label{eq11}\end{equation}

where we omitted for convenience the tilde accent on the new variables
$x_{i}$ in the transformed equation. This is exactly the Lotka-Volterra
form announced in (\ref{eq8}) with $M=BA$, i.e. M is the invariant
matrix (\ref{eq7}). \\

For the non-square QP systems, i.e. the systems whose matrices $A$
and $B$ respectively are $n\times N$ and $N\times n$, the Lotka-Volterra
form (\ref{eq8}) is also shown to be the canonical form with a matrix
$M=BA$ that is $N\times N$ \cite{2,6}. Thus, the reduction to the
Lotka-Volterra canonical form is a general property of the QP dynamical
systems.\\

The fact that the transformation leading to the canonical form is
a diffeomorphism reduces the analysis of the solutions of general
QP-systems to those of the corresponding Lotka-Volterra systems. This
result paves the way to the transfer of the whole corpus of knowledge
related to the latter equations which is pretty extensive.\\

We now show that this reduction to the Lotka-Volterra format leads
to an explicit formula for the general solution of the QP-systems.
This solution is expressed in terms of a Taylor series whose general
coefficient is explicitely known.

\section*{General solution of QP-systems}

The Taylor series in time $t$ for the solution $x_{i}(t)$ of the
Lotka-Volterra differential system (\ref{eq8}) is defined as 

\begin{equation}
x_{i}(t)=\sum_{k=0}^{\infty}c_{i}(k)\frac{t^{k}}{k!}\qquad\mbox{ for }i=1,\ldots,N\label{eq12}\end{equation}

The coefficient $c_{i}(k)$ is calculated by performing the $k$ order
time derivative of $x_{i}(t)$ at time $t=0$. This is readily done
by iterating the $t$ derivative while recursively using the Lotka-Volterra
system (\ref{eq8}). The result is amazingly simple: \begin{eqnarray}
c_{i}(k)=x_{i}(0)\sum_{i_{1}=0}^{N}\ldots\sum_{i_{k}=0}^{N}M_{ii_{1}}\left(M_{ii_{2}}+M_{i_{1}i_{2}}\right)\ldots\left(M_{ii_{k}}+M_{i_{1}i_{k}}+\cdot\cdot\cdot+M_{i_{k-1}i_{k}}\right)\nonumber \\ 
x_{i_{1}}(0)\ldots x_{i_{k}}(0)\label{eq13}\end{eqnarray}

where the $x_{i}(0)$ are the components of the initial condition.
Here again we omitted the tilde accent on the variables $x_{i}$ .\\

The considered Lotka-Volterra system results from the transformation
(\ref{eq8},\ref{eq9},\ref{eq10}) of a QP-system (\ref{eq2}).
Going back to the original variables of this QP-system, we get for
the order $k$ coefficient $C_{i}(k)$ of the Taylor series solving
this system: \begin{eqnarray}
C_{i}(k)=x_{i}(0)\sum_{i_{1}=0}^{N}\ldots\sum_{i_{k}=0}^{N}A_{ii_{1}}\left(A_{ii_{2}}+M_{i_{1}i_{2}}\right)\ldots\left(A_{ii_{k}}+M_{i_{1}i_{k}}+\ldots+M_{i_{k-1}i_{k}}\right)\nonumber \\
\prod_{j_{1}=1}^{n}x_{j_{1}}(0)^{B_{i_{1}j_{1}}}\ldots\prod_{j_{k}=1}^{n}x_{j_{k}}(0)^{B_{i_{k}j_{k}}}\label{eq14}\end{eqnarray}
 where $i$ runs from $1$ to $n$ and $n$ is the dimension of the
original QP-system. Thus, the general solution of the QP-system is
\begin{equation}
x_{i}(t)=\sum_{k=0}^{\infty}C_{i}(k)\frac{t^{k}}{k!}\qquad\mbox{ for }i=1,\ldots,n\label{eq15}\end{equation}

An observation of the coefficients $c_{i}(k)$ clearly reveals a combinatorial
structure related to the product \begin{equation}
M_{ii_{1}}\left(M_{ii_{2}}+M_{i_{1}i_{2}}\right)\dots\left(M_{ii_{k}}+M_{i_{1}i_{k}}+\ldots+M_{i_{k-1}i_{k}}\right)\label{eq16}\end{equation}
. 

For a one dimensional Lotka-Volterra system this product would reduce
to $k!M^{k}$. In more dimensions the number $M$ is replaced by the
components of the matrix $M_{ij}$. By doubling the number of sums
and introducing the corresponding Kronecker deltas the previous product
becomes \begin{equation}
\sum_{j_{1}=0}^{N}\dots\sum_{j_{k}=0}^{N}\delta_{ij_{1}}\left(\delta_{ij_{2}}+\delta_{i_{1}j_{2}}\right)\dots\left(\delta_{ij_{k}}+\delta_{i_{1}j_{k}}+\ldots+\delta_{i_{k-1}j_{k}}\right)M_{j_{1}i_{1}}M_{j_{2}i_{2}}\ldots M_{j_{k}i_{k}}\label{eq17}\end{equation}
 where the tensor \begin{equation}
\delta_{ij_{1}}(\delta_{ij_{2}}+\delta_{i_{1}j_{2}})\dots\left(\delta_{ij_{k}}+\delta_{i_{1}j_{k}}+\ldots+\delta_{i_{k-1}j_{k}}\right)\label{eq18}\end{equation}
 is a N-dimensional generalization of the factorial function $k!$.
Its role in the structure of the Taylor coefficient is double. It
governs the type of contracted products between the $k$ tensors $M$
that appear in the coefficient $c_{i}(k)$, and also counts these
types. The combinatorial study of this object would certainly lead
to interesting perspectives. \\

The above Taylor series (\ref{eq12}), when it converges, can be considered as defining
a new class of special functions which presents some analogy with the large
class of the hypergeometric functions. In analogy with the latter class of special functions, the study of other properties of  these new functions such as asymptotic series expansions and integral representations would be of the greatest interest. Indeed,
the QP-systems are known to have for certain values of their parameters
chaotic solutions. The latter have extremely complex behaviours in $t$ such
as the sensitivity with respect to the variations of the initial conditions or 
strange attractors having fractal geometry in
the phase-space. These properties should be related in some way to
the properties of the functions that we generated as solutions of
these systems. Such a study would certainly contribute to develop
a bridge between the theory of nonlinear dynamical systems and the theory of
special functions. \\
Furthermore, we believe that the study of the combinatorial properties of the tensor (\ref{eq18}) could be of great interest for both combinatorics and dynamical system theory.


\begin{thebibliography}{7}
\bibitem[1]{1}L.Brenig, Complete factorisation and analytic solutions
of generalized Lotka-Volterra equations, Phys.Lett. \textbf{A133},
378 (1988). 

\bibitem[2]{2}L.Brenig and A.Goriely, Universal canonical forms for
time-continuous dynamical systems, Phys.Rev. \textbf{A40}, 4119 (1989). 

\bibitem[3]{3}M.Peschel and W.Mende, The predator-prey model. Do
we live in a Volterra world?, Springer, New York, 1986. 

\bibitem[4]{4}J.L.Gouz\'e, Transformations of polynomial differential
systems in the positive octant, Rapport INRIA n\textdegree{}1308,
Sophia-Antipolis, France (1990).

\bibitem[5]{5}J.L.Gouz\'e, Global behaviour of n-dimensional Lotka-Volterra
systems, Mathematical Biosciences {\bf 113}, 231 (1993). 

\bibitem[6]{6}A.Goriely and L.Brenig, Algebraic degeneracy and partial
integrability for systems of ordinary differential equations, Phys.Lett.
{\bf A145}, 245 (1990). 

\bibitem[7]{7}E.H.Kerner, Universal formats for nonlinear ordinary
differential equations, J.Math.Phys. {\bf 22}, 1366 (1981).
\end{thebibliography}
\end{document}